\documentclass{amsart}
\usepackage{amssymb}
\input{secondary.sty}

\newtheorem{thm}{THEOREM}[section]
\newtheorem{conj}[thm]{CONJECTURE}
\newtheorem{cor}[thm]{COROLLARY}
\newtheorem{defn}[thm]{DEFINITION}
\newtheorem{ex}[thm]{EXAMPLE}

\newtheorem{prob}[thm]{PROBLEM}
\newtheorem{prop}[thm]{PROPOSITION}

\begin{document}

\title{Characteristic classes for Riemannian foliations}

\author{Steven Hurder}
\thanks{Supported in part by   NSF grant  DMS-0406254}
\address{Department of Mathematics, University of Illinois at Chicago,    851 S. Morgan St., Chicago, IL 60607-7045, USA}
\email{hurder@uic.edu}
\thanks{Preprint date: October 26, 2008; rev. December 8, 2008}

\date{}



\begin{abstract}
 The purpose of this paper is to both survey and offer some new results on the non-triviality of the characteristic classes of Riemannian foliations. We give  examples where the primary Pontrjagin classes are all linearly independent. The independence of the secondary classes is also discussed, along with their total variation. Finally, we give a negative solution of a conjecture that the map of classifying spaces $\FRGq \to \FGq$ is trivial for codimension $q > 1$.
\end{abstract}

\maketitle

\section{Introduction} \label{sec-intro}
The Chern-Simons class \cite{ChernSimons1974} of a closed $3$-manifold $M$, considered as foliated by its points, is the most well-known of the secondary classes   for Riemannian foliations.  Foliations with leaves of positive dimension offer a much richer class to study, and the values of their secondary classes reflect both geometric (metric) and dynamical properties of the foliations.   It is known that all of these classes can be realized independently for explicit examples (Theorem~\ref{thm-injectall}), but there remain a number of open problems to study. The purpose of this note is to survey the known results,   highlight some of the questions that remain unanswered, and also to provide a negative answer to one of these questions.

Let $M$ be a   smooth manifold of dimension $n$, and let $\F$ be  a smooth foliation  of     codimension $q$.    We say that $\F$ is a \emph{Riemannian foliation} if   there is a smooth  Riemannian metric $g$ on $TM$ which is \emph{projectable} with respect to $\F$. Identify the normal bundle  $Q$   with the orthogonal space $T\F^{\perp}$, and let $Q$ have the restricted Riemannian metric $g_Q =  g | Q$.  For a vector $X \in T_xM$ let $X^{\perp} \in Q_x$ denote its orthogonal  projection. Given a leafwise path $\gamma$ between points $x, y$ on a leaf $L$, the transverse holonomy $h_{\gamma}$ along $\gamma$ induces a   linear transformation $dh_x[\gamma] \, \colon Q_x \rightarrow Q_y$. The   fact that the Riemannian metric $g$ on $TM$ is projectable  is equivalent to the fact  that    the transverse  linear holonomy transformation $dh_x[\gamma]$ is an isometry for all such paths \cite{Haefliger1985, Haefliger1989, MM2003, Molino1982, Molino1988, Reinhart1983}.

There are a large variety of examples of Riemannian foliations which arise naturally in geometry. Given a smooth fibration $\pi \colon M \to B$, the connected components of the fibers of $\pi$ define the leaves of a foliation $\F$ of $M$. A Riemannian metric $g$ on $TM$ is projectable if there is a Riemannian metric $g_B$ on $TB$ such that the   restriction of $g$  to the normal bundle $Q \equiv T\F^{\perp}$  is the lift of the metric   $g_B$. The pair $(\pi \colon M \to B, g)$ is said to be a \emph{Riemannian submersion}. Such foliations provide  the most basic examples of   Riemannian foliations. 

 Suspensions of isometric actions of finitely generated groups provide another canonical class of examples of Riemannian foliations. The celebrated Molino Structure Theory for Riemannian foliations of compact manifolds reduces, in a broad sense,  the study of the geometry of Riemannian foliations to a m\'{e}lange of these two types of examples -- a combination of fibrations and group actions; see   Theorems~\ref{thm-molino1} and \ref{thm-molino2} below.  When the dimension of $M$ is at most $4$, the Molino approach yields  a ``classification'' of all Riemannian foliations. However, in general the structure theory is too rich and subtle to effect a classification for codimension $q \geq 3$ and leaf dimensions $p \geq 2$. The survey by Ghys, Appendix E of \cite{Molino1988}, gives an overview of the classification problem circa 1988.

The secondary characteristic classes of Riemannian foliations give another approach, independent of the Molino results,  for a broad  classification of Riemannian foliations.  While the secondary classes do not provide as precise a classification scheme as the Molino Theory, their study focuses  attention on various classes of Riemannian foliations, which can then be investigated   in terms of known examples, their  Molino Structure Theory, and the values of their characteristic classes,   often leading to new insights.

The characteristic classes of a Riemannian foliation are divided into three types: the primary classes, given by the ring generated by the  Euler and Pontrjagin classes; the secondary classes; and the blend of these two as defined by the  Cheeger-Simons differential characters.  Each of these types of invariants have  been more or less extensively studied, and part of the goal of this paper is to survey many of the known results. The paper also includes various new results and unpublished observations, some of which were presented in the author's talk \cite{Hurder2003a}. In a subsequent work    \cite{HurderToeben2008}, we will relate the values of the characteristic classes with the geometry of Riemannian foliations and the Molino classification.
 
The main new result  of this paper uses  characteristic classes to give a negative answer to Conjecture~3 of the Ghys survey [{\it op. cit.}]. The proof of the following is given in \S\ref{sec-primary}.

\begin{thm} 
For $q \geq 2$, the map $H_{4k-1}(\FRGq ; \mZ) \to H_{4k-1}(\FGq ; \mZ)$ has infinite-dimensional image for all   degrees  $4k -1 \geq q$.
\end{thm}

This paper is an expanded  version of a talk given at the joint AMS-RSME Meeting in Seville, Spain in June 2003. 
The talk was dedicated to the memory of Connor Lazarov, who   passed away on February 27, 2003. We dedicate this work to his memory, and especially  his fun-loving approach to all things, including his mathematics, which contributed so much to the field of Riemannian foliations.

\bigskip

 \section{Classifying spaces} \label{sec-classifying}

 The universal Riemannian groupoid $R\Gamma_q$ of $\mR^q$ is    the groupoid generated by the collection of all local isometries 
$\gamma \colon (U_{\gamma}, g'_{\gamma})  \to (V_{\gamma}, g''_{\gamma})$ where $g'_{\gamma}, g''_{\gamma}$ are complete Riemannian metrics on $\mR^q$, and $U_{\gamma}, V_{\gamma} \subset \mR^q$ are open subsets. The realization of the groupoid  $R\Gamma_q$ is a Hausdorff  topological space $BR\Gamma_q$, which is well-defined up to weak-homotopy equivalence \cite{Haefliger1970, Haefliger1971}. If we restrict to orientation-preserving maps of $\mR^q$, then we obtain the groupoid denoted by $R\Gamma^+_q$ with classifying space $BR\Gamma^+_q$.

 The universal   groupoid $\Gamma_q$ of $\mR^q$ is that defined by the groupoid generated by the collection of all local diffeomorphisms  
$\gamma \colon U_{\gamma}  \to V_{\gamma}$ where   $U_{\gamma}, V_{\gamma} \subset \mR^q$ are open subsets. The realization of the groupoid  $\Gamma_q$ is a non-Hausdorff  topological space $B\Gamma_q$, which is well-defined up to weak-homotopy equivalence.
 
An $R\Gamma_q$--structure on $M$ is an open covering ${\mathcal U} = \{U_{\alpha} \mid \alpha \in {\mathcal A} \}$ of $M$  
and for each $\alpha \in {\mathcal A}$, there is given  
\begin{itemize}
\item a smooth map $f_{\alpha} \colon U_{\alpha} \to V_{\alpha} \subset \mR^q$
\item a Riemannian metric $g'_{\alpha}$ on $\mR^q$ 
\end{itemize}
such that the pull-backs $f_{\alpha}^{-1}(T\mR^q) \to U_{\alpha}$ define a smooth vector bundle $Q \to M$ with Riemannian metric $g |  Q = g_{\alpha} = f_{\alpha}^*g'_{\alpha}$.   An $R\Gamma_q$--structure on $M$ determines a continuous map $M  \to |\cU | \to   BR\Gamma_q$.  

 Foliations $\F_0$ and $\F_1$ of codimension $q$ of $M$ are {\it integrably homotopic} if there is a foliation $\F$ of $M \times \bR$ of codimension $q$ such that $\F$ is everywhere transverse to the slices $M \times \{t\}$, so defines a foliation $\F_t$ of codimension-$q$ of $M \times \{t\}$, and $\F_t$ of $M_t$ agrees with $\F_t$ of $M$ for $t =0,1$.  This notion extends to Riemannian foliations, where we require that $\F$ defines a Riemannian foliation of codimension-$q$ of $M \times \mR$.
  
For the following,  assume that the normal bundle $Q \to M$ of $\F$ is orientable.

\begin{thm} [Haefliger \cite{Haefliger1970, Haefliger1971}] A Riemannian foliation $(\F, g)$ of $M$ with oriented normal bundle defines an $R\Gamma^+_q$--structure on $M$.  The homotopy class of the composition
$h_{\F, g} \colon M  \to BR\Gamma^+_q$
depends only on the integrable homotopy class of $(\F, g)$.
\end{thm}

The derivative of the local maps $\gamma \colon (U_{\gamma}, g'_{\gamma})  \to (V_{\gamma}, g''_{\gamma})$ takes values in $\SOq$, and induces a  classifying  map $\nu \colon BR\Gamma^+_q \to B\SOq$. The homotopy fiber of $\nu$ is denoted by    $\FRGq$.  The space  $\FRGq$  classifies $R\Gamma^+_q$--structures with a (homotopy class of) framing for $Q$. Let $P \to M$ be the bundle of oriented orthonormal frames of $Q \to M$, and $s \colon M \to P$ a choice of framing of $Q$. Then we have the commutative diagram:
\bigskip

\begin{picture}(500,150)\label{eqn-diag}
\put(95,140){${\SOq}$}
\put(180,140){${\SOq}$}
\put(290,140){$\SOq$}

\put(145,140){$=$}
\put(255,140){$=$}

\put(105,120){$\downarrow$}
\put(195,120){$\downarrow$}
\put(300,120){$\downarrow$}

\put(140,110){$\overline{h^s_{\F,g}}$}

\put(100,100){$P$}
\put(180,100){$\FRGq$}
\put(295,100){$\FGq$}

\put(140,100){$\longrightarrow$}
 \put(250,100){$\longrightarrow$}

\put(90,80){$s$}
\put(100,80){$\uparrow$}
\put(110,80){$\downarrow$}
\put(195,80){$\downarrow$}
\put(300,80){$\downarrow$}

\put(140,65){$h_{\F, g}$}
\put(100,55){$M$}
\put(180,55){$BR\Gamma^+_q$}
\put(295,55){$B\Gamma^+_q$}

\put(145,55){$\longrightarrow$}
\put(255,55){$\longrightarrow$}

\put(195,35){$\downarrow$}
\put(215,35){$\nu$}
\put(300,35){$\downarrow$}
\put(310,35){$\nu$}

\put(180,15){$B\SOq$}
\put(290,15){$B\SOq$}
\put(255,15){$=$}
\end{picture}

where the right-hand column is the corresponding sequence of classifying spaces for the groupoid defined by the germs of  local diffeomorphisms of $\mR^q$. The natural  maps $\FRGq \to \FGq$ and  $BR\Gamma^+_q \to B\Gamma^+_q$ are induced by the natural transformation which ``forgets'' the normal Riemannian metric data.

The approach to classifying foliations initiated by Haefliger in \cite{Haefliger1970,Haefliger1971}  is based on the study of the homotopy classes of maps $[M, BR\Gamma^+_q]$ from a given manifold $M$ without boundary to $BR\Gamma^+_q$. Given   a   homotopy class of an embedding of an oriented subbundle $Q \subset TM$ of dimension $q$, one studies the homotopy classes of  maps  $h_{\F, g} \colon M \to BR\Gamma^+_q$ such  that  the composition  $\nu \circ h_{\F,g} \colon M \to B\SOq$     classifies the homotopy type of the   subbundle $Q$.  This result motivates the study of the homotopy types of the spaces $BR\Gamma^+_q$ and $\FRGq$.

In the case of codimension-one, a  Riemannian foliation with oriented normal bundle of $M$ is equivalent to specifying a closed, non-vanishing 1-form $\omega$  on $M$. As $\mathbf{SO}(1)$ is the trivial group, $FR\Gamma^+_1 = BR\Gamma^+_1$, and the classifying map $M \to BR\Gamma^+_1$ is determined by the real cohomology class of $\omega$, which follows from the following result of Joel Pasternack. Let $\mR_{\delta}$ denote the real line, considered as a \emph{discrete} group, and $B{\mR}_{\delta}$ its classifying space.

\begin{thm}[Pasternack \cite{Pasternack1973}] \label{thm-pasternack1}   
 There is a natural homotopy equivalence  $BR\Gamma^+_1  \simeq FR\Gamma_1  \simeq B{\mR}_{\delta}$.
\end{thm}

For codimension $q \geq 2$, the space $B\SOq$ is no longer contractible, as $H^*(B\SOq ; \mZ)$ is generated as a commutative algebra generated by the Pontrjagin classes, and for $q$ even, by the Pontrjagin classes along with the Euler class in degree $q$. This will be discussed further in the next section. 

Pasternack's Theorem admits the following partial generalization:

\begin{thm}[Hurder \cite{Hurder1980,Hurder1981a}] \label{thm-connected}   Let $q \geq 2$. Then the space $\FRGq$ is $(q-1)$--connected. That is, $\pi_{\ell}(\FRGq) = \{0\}$ for $0 \leq \ell < q$. Moreover, the volume form associated to the transverse metric defines a surjection $\vol \colon \pi_q(\FRGq) \to \mR$.
 \end{thm}
\sop We just give a sketch; details appeared in \cite{Hurder1981a}. 
Following a remark by Milnor, one observes that by the Phillips Immersion Theorem \cite{Phillips1967, Phillips1968, Phillips1969}, an $FR\Gamma_q$--structure on ${\mS}^{\ell}$ for $0< \ell < q$ corresponds to 
a Riemannian metric defined on an open neighborhood retract of the $\ell$-sphere,  ${\mS}^{\ell} \subset U \subset {\mR}^q$. 

 Given an $FR\Gamma_q$--structure on  the open set  $U \subset {\mR}^q$ -- which is equivalent to specifying a Riemannian metric on $TU$ -- 
 one then constructs an explicit integrable homotopy through framed $R\Gamma_q$-structures on a smaller open neighborhood  ${\mS}^{\ell} \subset V \subset U$.  The integrable homotopy starts with the given Riemannian metric on $TV$, and ends with  the standard Euclidean metric on   $TV$, which represents the  ``trivial'' $FR\Gamma_q$--structure on  ${\mS}^{\ell}$.    Thus, every $FR\Gamma_q$--structure on ${\mS}^{\ell}$ is homotopic to the trivial structure.

The surjection $\vol \colon \pi_q(\FRGq) \to \mR$ is well-known, and is realized simply by varying the total volume of a Riemannian metric on $\mS^q$, foliated by points.  \hfill $\eop$

Associated to the classifying map $\nu \colon BR\Gamma^+_q  \to  B\SOq$ is the Puppe sequence
\begin{equation}\label{eq-puppe}
 \cdots \longrightarrow \Omega  \FRGq   \stackrel{\Omega{\nu}}{\longrightarrow} \Omega  BR\Gamma^+_q  \longrightarrow \SOq \stackrel{\delta}{\longrightarrow}  \FRGq \longrightarrow   BR\Gamma^+_q   \stackrel{\nu}{\longrightarrow}  B\SOq
\end{equation}
In the case of codimension $q = 2$, ${\mathbf{SO}(2)} = \mS^1$ and $FR\Gamma_2$ is $1$-connected, so the map 
$\delta \colon {\mathbf{SO}(2)} \to FR\Gamma_2$ is contractible. This yields as an immediate consequence:
 
\begin{thm}[Hurder \cite{Hurder1993}] \label{thm-puppe}   
$\Omega BR\Gamma^+_2 \cong {\mathbf{SO}(2)} \times \Omega FR\Gamma_2$. 
 \end{thm}
It is noted in \cite{Hurder1993}  that the homotopy equivalence in Theorem~\ref{thm-puppe} is not an $H$-space equivalence, as this would imply that map  
$\nu^* \colon H^*(B{\mathbf{SO}(2)} ; \mR) \to  H^*(BR\Gamma_2 ; \mR)$ is an injection, which is false.
In contrast, we have the following result: 

\begin{thm} \label{thm-loops} The connecting map 
$\delta \colon  \SOq \to \FRGq$ in (\ref{eq-puppe})  is not homotopic to a constant for   $q \geq 3$. 
\end{thm}
 Note that the map $\delta \colon  \SOq \to \FRGq$
classifies the   Riemannian foliation with standard framed normal bundle on $\SOq \times \mR^q$, obtained via the pull-back of the standard product foliation of $\SOq \times \mR^q$  via the action  of    $\SOq$ on $\mR^q$.  Theorem~\ref{thm-loops}  asserts that the canonical twisted   foliation of $\SOq \times \mR^q$   is  not integrably homotopic  through framed Riemannian foliations  to the standard product foliation. This will be proven in section~\ref{sec-secondary}, using basic properties of the secondary classes for Riemannian foliations. For the non-Riemannian case, it is conjectured that the connecting map $\delta \colon  \SOq \to \FGq$ is homotopic to the constant map \cite{Hurder2003b}.

 To close this discussion of general properties of    the classifying spaces of Riemannian foliations, we pose a problem particular to codimension two:   

\begin{prob} \label{prob-two}
Prove that the map induced by the volume form
 $\vol \colon \pi_2(FR\Gamma_2) \to \mR$ is an isomorphism. That is, given two $R\Gamma_2$-structures $\F_0$ and $\F_1$ on $M = \mR^3 - \{0\}$, with homotopic normal bundles, prove that $\F_0$ and $\F_1$ are   homotopic as  $R\Gamma_2$-structures if and only if they have cohomologous transverse volume forms.
 \end{prob}
One can view this as asking for a ``transverse uniformization theorem'' for Riemannian foliations of codimension two. 
Note that Example~\ref{ex-LP2} below shows the conclusion of Problem~\ref{prob-two}  is false for $q=3$.

 \section{Primary classes} \label{sec-primary}
 
 The primary classes of a  Riemannian foliation are those obtained from the cohomology of the classifying space of the normal bundle $Q \to M$, pulled-back via the classifying map $\nu \colon M \to B\SOq$. Recall \cite{MilnorStasheff1974} that the cohomology groups of $\SOq$ are isomorphic to   free polynomial ring:
\begin{eqnarray*}
H^*(B{\mathbf{SO}(2)} ; \mZ) & \cong & \mZ[E_1]\\
H^*(B\SOq ; \mZ) & \cong & \mZ[E_m, P_1, \ldots , P_{m-1}] ~, ~ q = 2m \geq 4\\
H^*(B\SOq ; \mZ) & \cong & \mZ[P_1, \ldots , P_m] ~, ~ q = 2m + 1 \geq 3
\end{eqnarray*}
As usual, $P_{j}$ denotes the   Pontrjagin cohomology class of degree $4 j$, $E_m$ denotes the   Euler class of degree $2m$, and the square $E_m^2 = P_m$ is the top degree generator of the Pontrjagin ring.
   
  There are three main results concerning the universal map  $\nu^* \colon H^{\ell}(B\SOq ; {\mathcal R}) \to H^{\ell} (BR\Gamma^+_q ; {\mathcal R})$, where ${\mathcal R}$ is a coefficient ring, which we discuss in detail below.
  
  \begin{thm} [Pasternack \cite{Pasternack1973}] \label{thm-pasternack2}
  $\nu^* \colon H^{\ell}(B\SOq ; \mR) \to H^{\ell} (BR\Gamma^+_q ; \mR)$ is trivial for $\ell > q$.
  \end{thm}
  
  \begin{thm} [Bott, Heitsch \cite{BottHeitsch1972}] \label{thm-BH}
  $\nu^* \colon H^{*}(B\SOq ; \mZ) \to H^{*} (BR\Gamma^+_q ; \mZ)$ is injective.
  \end{thm}
  
  \begin{thm} [Hurder \cite{Hurder1981a,Hurder2003a}] \label{thm-thom}
  $\nu^* \colon H^{\ell}(B\SOq ; \mR) \to H^{\ell} (BR\Gamma^+_q ; \mR)$ is injective for $\ell \leq q$.
  \end{thm}
  
  The contrast between Theorems~\ref{thm-pasternack2} and \ref{thm-BH} is one of the themes of this section, while the proof of Theorem~\ref{thm-thom} is   based on an observation. 
   
Let  $\nabla_g$ denote the   Levi-Civita  connection on $Q \to M$ associated to  the projectable metric $g$ for $\F$. The Chern-Weil construction associates to each universal class $P_{j}$ the closed  Pontrjagin form  $p_{j}(\nabla_g) \in \Omega^{4 j}(M ; \mR)$.  For $q=2m$, as $Q$ is assumed to be oriented, there is also the Euler form $e_m(\nabla_g) \in \Omega^{2m}(M ; \mR)$ whose square $e_m(\nabla_g)^2 = p_m(\nabla_g)$. By Chern-Weil Theory,    the universal map 
$\nu^* \colon H^{*}(B\SOq ; \mR) \to H^{*} (BR\Gamma^+_q ; \mR)$ is defined by setting
$\nu^*(P_{j}) = [p_{j}(\nabla_g)] \in H^{4j}(M ; \mR)$, where $[ \beta ]$ represents the de~Rham cohomology class of a closed form $\beta$.

Let $m$ be the least integer such that $q \leq 2m+2$. Given   $J= (j_1,  j_2,  \ldots , j_m)$ with each $j_{\ell} \geq 0$, set 
    $p_J = p_1^{j_1} \cdot p_2^{j_2} \cdots p_m^{j_m}$, which has degree    $4|J| = 4(j_1+ \cdots + j_m)$. Then let $\cP$ denote a general monomial, which for $q = 2m+1$ has the form $\cP = p_J$ with $\deg(\cP) = 4 |J|$. 
For  $q = 2m$,  either let $\cP = p_J$ with  $\deg(\cP) = 4 |J|$, or $\cP = e_m \cdot p_J$ with $\deg(\cP) = 4|J| + 2m$.

Pasternack   \cite{Pasternack1970} first observed  that the proof of the Bott Vanishing Theorem \cite{Bott1970} can be strengthened in the case of Riemannian foliations, as the adapted metric $\nabla_g$ is projectable. He showed that on the level of differential forms,

  \begin{thm} [Pasternack \cite{Pasternack1970,Pasternack1973}] \label{thm-pasternack3} If $\deg(\cP) > q$ then 
$\cP(\nabla_g)  = 0$. 
\end{thm}
 Theorem~\ref{thm-pasternack2} follows immediately. Today, this result is considered ``obvious'', but that is due to the later extensive development of this field in the 1970's. 

Let us next consider the injectivity of the map  $\nu^* \colon H^{k}(B\SOq ; \mR) \to H^{k} (BR\Gamma^+_q ; \mR)$. We  recall a basic observation  of Thom: 

  \begin{thm} [\cite{MilnorStasheff1974}] \label{thm-thom2} 
  There is a compact, orientable  Riemannian manifold  $B$ of dimension $q$ such that all of the Pontrjagin and Euler classes up to degree $q$ are independent in $H^*(B ; \mR)$. 
If $q$ is odd, then $B$ can be chosen to be a connected manifold.
\end{thm}
\sop
For $q$ even, let $B$ equal the disjoint union of all products of the form 
$${\bf CP}^{i_1} \times \cdots \times {\bf CP}^{i_k} \times S^1 \times \cdots \times S^1$$
with dimension $q$.
For $q$ odd,   $B$ is the connected sum   of all products of the form 
$${\bf CP}^{i_1} \times \cdots \times {\bf CP}^{i_k} \times S^1 \times \cdots \times S^1$$
with dimension $q$. The claim then follows by the Splitting Principle \cite{MilnorStasheff1974}  for the Pontrjagin classes. \hfill $\eop$

The claim of Theorem~\ref{thm-thom} now follows from the universal properties of $BR\Gamma^+_q$, as we endow the manifold $B$ with the foliation $\F$ by points, with  the standard   Riemannian metric on $B$. \hfill $\eop$

The original proof of Theorem~\ref{thm-thom} in \cite{Hurder1981a} used the fact that $\nu \colon BR\Gamma^+_q \to B\SOq$ is $q$-connected.

\medskip
 
Next, we discuss the results of Bott and Heitsch from \cite{BottHeitsch1972}. Let $K \subset \SOq$ be a closed Lie subgroup, and let $\Gamma \subset K$ be a finitely-generated subgroup.  Suppose that $B$ is a closed connected manifold, with basepoint $b_0 \in B$,  for which there is a surjection of its fundamental group $\rho \colon \Lambda = \pi_1(B, b_0) \to \Gamma \subset K \subset \SOq$. Then via the natural action of $\SOq$ on $\mR^q$ we obtain an action of $\Lambda$ on $\mR^q$. Let $\widetilde{B} \to B$ denote the universal covering of $B$, equipped with the right action of $\Lambda$ by deck transformations.  Then form the flat bundle 
\begin{equation}
\mE_{\rho} = \widetilde{B} \times \mR^q / (b \cdot \gamma, \vec{v}) \sim (b, \rho(\gamma) \cdot \vec{v}) \stackrel{\pi}{\longrightarrow} \widetilde{B}/\Lambda = B
\end{equation}
As the action of $\Lambda$ on $\mR^q$ preserves the standard Riemannian metric, we obtain a Riemannian foliation $\F_{\rho}$ on $\mE_{\rho}$ whose leaves are the integral manifolds of the flat structure. The classifying map of the foliation $\F_{\rho}$ is given by the composition of maps
\begin{equation}
\mE_{\rho} \to B\Lambda \to B(K_{\delta}) \to B(\SOq_{\delta}) \to BR\Gamma^+_q
\end{equation}
where $K_{\delta}$ and $\SOq_{\delta}$ denotes the corresponding Lie groups   considered with the discrete topology, and $B(K_{\delta})$ and $B(\SOq_{\delta})$ are the corresponding classifying spaces. 

The Bott-Heitsch examples take $K$ to be a maximal torus, so that for $q =  2m$ or $q = 2m+1$, we have $K = \mT^m =  \SOtwo \times \cdots \times \SOtwo$ with $m$ factors. Consider first the case $q = 2$. For an odd prime $p$, let  $\Gamma = \mZ/p \mZ$,  embedded as the $p$-th roots of unity in $K = \SOtwo$. Let $B = \mS^{2\ell + 1}/\Gamma$ be the quotient of the standard odd-dimensional sphere, and consider the composition
\begin{equation}
\nu \circ \rho \colon B \to \mE_{\rho} \to B(\mZ/p \mZ) \to     B(\SOtwo_{\delta}) \to BR\Gamma^+_2 \stackrel{\nu}{\longrightarrow} B\SOtwo
\end{equation}
The composition $\nu \circ \rho$ classifies the Euler class of the flat bundle $\mE_{\rho} \to B$, which is torsion.  The map in cohomology with $\mZ/ p \mZ$-coefficients,
\begin{equation}
(\nu \circ \rho)^* \colon H^{*}(B\SOtwo ; \mZ/ p \mZ) \to H^{*}(B; \mZ/p \mZ)
\end{equation}
  is injective for $* \leq  2 \ell$. It follows that the map 
\begin{equation}
\nu^* \colon H^{*}(B\SOtwo ; \mZ/ p \mZ) \to H^{*}(BR\Gamma^+_2; \mZ/p \mZ)
\end{equation}
is injective in all degrees. As this holds true for all odd primes,   it is also  injective for integral cohomology. \hfill $\Box$

Theorem~\ref{thm-BH} is a striking result, as     Theorem~\ref{thm-pasternack2} states that 
$\nu^* \colon H^{*}(B\SOtwo ; \mQ) \to H^{*}(BR\Gamma^+_2; \mQ)$ is the trivial map for $* > 2$. One thus concludes from  the Universal Coefficient Theorem for cohomology, that the  homology groups $H_*(BR\Gamma^+_2 ; \mZ)$ cannot be finitely generated in all odd degrees $* \geq 3$ (cf. \cite{BottHeitsch1972}).

The treatment of the cases where $q = 2m > 2$ or $q = 2m+1 > 2$ follows similarly, where one takes $\Gamma = (\mZ/ p \mZ)^m \subset \mT^m \subset \SOq$, and let $p \to \infty$. An application of the splitting theorem for the Pontrjagin classes of vector bundles then yields Theorem~\ref{thm-BH}. 

The fibration sequence
$\ds  \FRGq \to BR\Gamma^+_q \to B\SOq$ 
yields a spectral sequence converging to the homology groups $ H^{*}(BR\Gamma^+_q; \mZ)$ with 
$E^2$-term 
$$E^2_{r,s} \cong H_r(B\SOq ; H_s(\FRGq ; \mZ))$$
It follows that the groups $H_s(\FRGq ; \mZ)$ cannot all be finitely generated for odd degrees $* \geq q$. In fact, we will see that this follows from the results of Pasternack and Lazarov discussed in the next section on secondary classes, but the homology classes being detected via the torsion classes above   seem to be of a different ``sort'' than those detected via the secondary classes. 

Recall that the universal  classifying map $\FRGq \to \FGq$ simply  ``forgets'' the added structure of a holonomy-invariant transverse Riemannian metric for the foliation. The following conjecture (page 308, \cite{Molino1988}) is that this map must be ``trivial'': 
\begin{conj} \label{conj-trivial}
For all $k > 0$, the      map of homotopy groups  $\pi_{k}(\FRGq) \to \pi_{k}(\FGq)$  is trivial. 
\end{conj}
The ideas of the  proof of Theorem~\ref{thm-BH} can be used to show this conjecture  is false.

\begin{thm}\label{thm-nontrivial}
For $2k > q \geq 2$, the image of $H_{4k-1}(\FRGq ; \mZ) \to H_{4k-1}(\FGq ; \mZ)$ is infinite-dimensional.
\end{thm}
\sop
Our approach uses the methods of the homological proof of Theorem~\ref{thm-BH} in \cite{BottHeitsch1972}, and especially    the   commutative  diagram from page 144 and the associated arguments. 

Let $\cP \in H^{4k}(B\SOq; \mZ)$ for  $4k > q$ be a generating monomial. The Bott-Heitsch Theorem~\ref{thm-BH} implies that the image $f^* \circ \nu^*(\cP) \in H^{4k}(BR\Gamma^+_q ; \mZ)$ is not a torsion class under the composition
\begin{equation}
H^{4k}(B\SOq; \mZ) \stackrel{\nu^*}{\longrightarrow}  
H^{4k}(B\Gamma^+_q ; \mZ)  \stackrel{f^*}{\longrightarrow}   
H^{4k}(BR\Gamma^+_q ; \mZ)
\end{equation}
where $f \colon BR\Gamma^+_q \to B\Gamma^+_q$ is the    ``forgetful'' map, forgetting the transverse Riemannian metric structure.

Let $\cA_{4k-1} = {\rm image} \{ H_{4k-1}(BR\Gamma^+_q ; \mZ) \to H_{4k-1}(B\Gamma^+_q ; \mZ) \}$.
Suppose that   $\cA_{4k-1}$ is finite-dimensional, then  ${\rm Ext} (\cA_{4k-1},  \mZ)$ is a torsion group.
Consider the commutative diagram:

   \medskip

\begin{picture}(420,160)

\put(170,140){$H^{4k}(B\SOq; \mZ)$}

\put(200,130){\vector(0,-1){16}}
\put(185,120){$\nu^*$}

\put(30,90){${\rm Ext} (H_{4k-1}(B\Gamma^+_q ; \mZ), ~ \mZ) $}
\put(175,90){$ H^{4k}(B\Gamma^+_q ; \mZ)$}
\put(270,90){${\rm Hom}(H_{4k}(B\Gamma^+_q ; \mZ), ~ \mZ)$}
\put(410,90){$\{0\}$}

\put(150,93){\vector(1,0){16}}
\put(240,93){\vector(1,0){16}}
\put(380,93){\vector(1,0){16}}

\put(154,98){$\tau$}
\put(246,98){$e$}

\put(70,80){\vector(0,-1){16}}
\put(200,80){\vector(0,-1){50}}
\put(320,80){\vector(0,-1){50}}

\put(50,70){$\iota^*$}
\put(50,30){$\sigma^*$}

\put(185,50){$f^*$}
\put(305,50){$f^*$}

\put(30,50){${\rm Ext} (\cA_{4k-1}, ~ \mZ) $}

\put(70,40){\vector(0,-1){16}}

\put(30,10){${\rm Ext} (H_{4k-1}(BR\Gamma^+_q ; \mZ), ~ \mZ) $}
\put(170,10){$ H^{4k}(BR\Gamma^+_q ; \mZ)$}
\put(265,10){${\rm Hom}(H_{4k}(BR\Gamma^+_q ; \mZ), ~ \mZ)$}
\put(410,10){$\{0\}$}

\put(150,13){\vector(1,0){16}}
\put(240,13){\vector(1,0){16}}
\put(380,13){\vector(1,0){16}}

\put(154,20){$\tau$}
\put(246,20){$e$}

\end{picture}

In the diagram, $e$ is the evaluation map of cohomology on homology, and $\tau$
maps onto its kernel. 

Also,  $\iota^*$ is induced by the inclusion $\iota \colon \cA_{4k-1} \subset  H_{4k-1}(B\Gamma^+_q ; \mZ)$, 
and $\sigma^*$ is induced by the surjection  $\sigma \colon  H_{4k-1}(BR\Gamma^+_q ; \mZ) \to  \cA_{4k-1}$.

\medskip

 The   Bott Vanishing Theorem  implies that the class 
 $$e \circ \nu^*(\cP) \in {\rm Hom}(H_{4k}(B\Gamma^+_q ; \mZ), ~ \mZ) \subset {\rm Hom}(H_{4k}(B\Gamma^+_q ; \mZ), ~ \mQ)$$
 is trivial for $\deg(\cP) > 2q$. Thus, there exists $\cP_{\tau} \in {\rm Ext} (H_{4k-1}(B\Gamma^+_q ; \mZ), ~ \mZ)$ such that 
 $\tau(\cP_{\tau}) = \nu^*(\cP)$.
 
 The class $\iota^*(\cP_{\tau}) \in {\rm Ext} (\cA_{4k-1}, ~ \mZ)$ is torsion, by the assumption on $\cA_{4k-1}$. Thus, 
 $f^* \circ \nu^*(\cP) = \tau \circ \sigma^* \circ \iota^* (\cP_{\tau})$ is a torsion class, which contradicts the Bott-Heitsch results. 
 
Thus, $\cA_{4k-1}$ cannot be finite-dimensional for  $4k > 2q$. \hfill $\eop$

\bigskip

 \begin{prob}
Find  geometric interpretations of the cycles in the image of the map 
 $$f_* \colon H_{4k-1}(\FRGq ; \mZ) \to H_{4k-1}(\FGq ; \mZ)$$
 \end{prob}
 The construction of foliations with solenoidal minimal sets in \cite{CH2008,Hurder2008c} give one realization of some of the classes in the image of this map, as discussed in the talk by the author \cite{Hurder2008b} at the conference of these Proceedings.  
 Neither the examples in \cite{Hurder2008c}, nor the situation overall, is   understood in sufficient depth.

  \bigskip

 \section{Secondary classes} \label{sec-secondary}
 
 Assume that $(\F , g)$ is a Riemannian foliation of codimension $q$. We also assume that  there exists a  framing of the normal bundle, denoted by $s \colon M \to P$ in section~\ref{sec-classifying}. 
 Then the data $(\F, g, s)$ yields a classifying map
 $h^s_{\F, g} \colon M \to \FRGq$. In this section, we discuss the construction of the secondary characteristic classes of such foliations, constructed using the Chern-Weil method as in \cite{Chern1979}, and some of the results about these classes. 

Recall that $\nabla_g$ denotes the Levi-Civita connection of the projectable metric on $Q$. 

Let $\cI(\SOq)$ denote the ring of Ad-invariant polynomials on the Lie algebra $\mathfrak{so}(q)$ of $\SOq$. Then we have 
\begin{eqnarray*}
\cI(\SOtwo)& \cong & \mR [e_m]\\
\cI(\mathbf{SO}(2m))& \cong & \mR [e_m, p_1, \ldots , p_{m-1}] ~, ~ q = 2m \geq 4\\
\cI(\mathbf{SO}(2m+1)) & \cong & \mR [p_1, \ldots , p_m] ~, ~ q = 2m + 1 \geq 3
\end{eqnarray*}
where the $p_{j}$ are the Pontrjagin polynomials, and $e_m$ is the Euler polynomial defined for $q$ even.

\vfill
\eject

The primary classes of $Q \to M$ are then represented by   closed forms,  obtained by applying the symmetric polynomials $p_j$ to the curvature matrix of 2-forms associated to the connection $\nabla_g$: 
$\Delta_{\F,g}(p_{j}) = p_{\ell}(\nabla_g) \in \Omega^{4j}(M)$. Set  $\Delta_{\F}[p_{j}] = [p_{j}(\nabla_g)] \in H^{4j}(M ; \mR)$. The Euler form $\Delta_{\F,g}(e_m) = e_m(\nabla_g) \in \Omega^{2m}(M)$ and the Euler class $\Delta_{\F,g}[e_m] \in H^{2m}(M, \mR)$ are similarly defined when $q = 2m$. We thus obtain   a multiplicative homomorphism 
$$\Delta_{\F,g} \colon \cI(\SOq) \to H^*(M ; \mR)$$

As noted in Theorem~\ref{thm-pasternack3}, Pasternack first observed that for $\nabla_g$ the adapted connection to a Riemannian foliation, the map $\Delta_{\F,g}$ vanishes identically in degrees greater than $q$.

\begin{defn}
For  $q=2m$,    set \quad  
$ \ds \cI(\SOq)_{2m} \equiv   \mR[e_m ,  p_1,p_2, \ldots, p_{m-1}]/ ( \cP \mid \deg (\cP) > q )$.

For   $q=2m+1$, set \quad
$\ds  \cI(\SOq)_{2m+1}  \equiv    \mR[p_1,p_2, \ldots, p_{m}]/  (\cP \mid \deg (\cP) > q )$. 
\end{defn}
 
 \medskip

\begin{cor} [Pasternack] Let $(\F,g)$ be a codiemsnion-$q$, Riemannian foliation of $M$. 
Then there is a well-defined characteristic homomorphism  $\Delta_{\F,g} \colon  \cI(\SOq)_q \to H^*(M ; \mR)$, which depends only upon the integrable homotopy class of $\F$. 
\end{cor}

 Of course, if we assume that the normal bundle $Q$ is trivial, then this map is zero in cohomology. The point of the construction of secondary classes is to obtain geometric information from the forms $p_{j}(\nabla_g) \in \Omega^{4j}(M)$ even if they are exact. If we do not assume that $Q$ is trivial, then one still knows that the cohomology classes $[p_{j}(\nabla_g)] \in H^{4j}(M ; \mR)$ lie in the image of the integral cohomology, $H^*(M ; \mZ) \to H^*(M ; \mR)$ so that one can use the construction of Cheeger-Simons differential characters as in \cite{CheegerSimons1973, ChernSimons1974, Koszul1975, Simons1972} to define secondary invariants in the groups $H^{4j -1}(M ; \mR/\mZ)$. These classes are closely related to the Bott-Heitsch examples above, and to the secondary classes constructed below.

Assume that we are given a trivialization 
 $s \colon M \to P$, then let 
  $\nabla_s$ be the   flat connection on $Q$ for which $s$ is parallel. 
Set $\nabla_t = t \nabla_g + (1-t) \nabla_s$, which we consider as a connection on the bundle $Q$ extended as product over $M \times \mR$. Then the Pontrjagin forms for $\nabla_t$ yield closed forms $p_{j}(\nabla_t) \in \Omega^{4j}(M \times \mR)$. Define the $4j -1$ degree transgression form
 \begin{equation}
 h_{j} = h_{j} (\nabla_g,s) = \int_0^1 \, \iota(\partial /\partial t) p_{j}(\nabla_t) \wedge dt ~  \in \Omega^{4j-1}(M)
\end{equation}
which satisfies the coboundary relation on forms:
$$ d h_{j}(\nabla_g,s) = p_{j} (\nabla_g) - p_{j}(\nabla_s) = p_{j} (\nabla_g) $$

For $q = 2m$ we also introduce the transgression of the Euler form, 
 \begin{equation}
\chi_m = \chi_m (\nabla_g,s) = \int_0^1 \, \iota(\partial /\partial t) e_m(\nabla_t) \wedge dt ~  \in \Omega^{q-1}(M) ~; ~ d \chi_m = e_m(\nabla_g)
\end{equation}

Note that if $4j > q$, then the form $p_{j}(\nabla_g) = 0$, so the transgression form $h_{j}$ is closed, and defines a secondary cohomology class $\Delta^s_{\F, g}(h_{j}) \in H^{4j -1}(M ; \mR)$.   In general, one follows the idea of the construction of the secondary classes for foliations to obtain the most general construction of invariants. Introduce the graded differential complexes, according to whether $q = 2m$ or $q=2m+1$:

$$
\begin{array}{lclcl}
RW_{2m}     &   = &    \Lambda \left(h_1, \ldots, h_{m-1}, \chi_m \right)    \otimes    \cI(\SOq)_{2m} &  , &   d_W(h_{j} \otimes 1) =  1 \otimes p_{j}   ~ , ~  d_W(  \chi_m \otimes 1 ) = e_m  \otimes 1 \\ 
RW_{2m+1}  & = &     \Lambda \left(h_1, \ldots, h_{m}  \right)    \otimes    \cI(\SOq)_{2m+1} &  , &      d_W(h_{j} \otimes 1) =  1 \otimes p_{j}  
 \end{array}
 $$
 
 For $I = (i_1 < \cdots < i_\ell)$ and $J= (j_1 \leq \cdots \leq j_k)$   set 
 \begin{equation}
h_I \otimes p_J =  h_{1_1} \wedge \cdots \wedge h_{i_{\ell}} \otimes p_{j_1} \wedge \cdots \wedge p_{j_k}
\end{equation}
Note that  $ \deg (h_I \otimes p_J) = 4(|I| + |J|) - \ell$, and that  $d_W( h_I \otimes p_J) = 0$ exactly when $4i_1 + 4 |J|  > q$. In the following, the expression $h_I \otimes p_J$ will always assume that the indexing sets $I$ and $J$ are ordered as above.

 \begin{thm}[Lazarov - Pasternack \cite{LazarovPasternack1976a}] 
 Let $(\F , g)$ be a Riemannian foliation of codimension $q \geq 2$ of a manifold $M$ without boundary, and assume that   there is given a framing of the normal bundle, 
  $s \colon M \to P$. Then the above constructions yield a map of differential graded algebras
  \begin{equation}
\Delta_{\F,g}^s \colon RW_q   \to   \Omega^*(M) 
\end{equation}
such that the induced map on cohomology, $\Delta_{\F,g}^s \colon H^*(RW_q )  \to   H^*(M ; \mR)$ is independent of the choice of basic connection $\nabla_g$, and depends only on the integrable homotopy class of $\F$ as a Riemannian foliation and the homotopy class of the framing $s$.
 \end{thm}
 
 As remarked by Kamber and Tondeur in \cite{KT1974d, KT1975a}, this construction can also be recovered from the method of truncated Weil algebras as applied to the Lie algebra $\mathfrak{so}(q)$. The functoriality of the construction of $\Delta_{\F,g}^s$ implies, in the usual way \cite{Lawson1975, LazarovPasternack1976a}:
 
 \begin{cor}
 There exists a universal characteristic homomorphism 
 \begin{equation}\label{eq-universal}
\Delta \colon H^*(RW_q, d_W) \to H^*(\FRGq; \mR)
\end{equation}
 \end{cor}

 There are many natural questions about how the values of these secondary classes are related to the geometry and dynamical properties of the foliation $(\F, g, s)$. We discuss some known results in the following.
 
 First, consider the role of the section $s \colon M \to P$. Given any smooth map $\varphi \colon M \to \SOq$, we obtain a new framing 
 $s' = s \cdot \varphi \colon M \to P$  by setting $s'(x) = s(x) \cdot  \varphi(x)$.
 Thus, $\varphi$   can be thought of as a gauge transformation of the normal bundle $Q \to M$.   
 
 The cohomology of the Lie algebra $\mathfrak{so}(q)$ is isomorphic to an exterior algebra, generated by the cohomology classes of left-invariant closed forms $\tau_{j} \in \Lambda^{4j -1}(\mathfrak{so}(q))$ for $j < q/2$, and the Euler form $\chi_m \in \Lambda^{2m -1}(\mathfrak{so}(q))$ when $q = 2m$. The   map $\varphi$ pulls these back to closed forms $\varphi^*(\tau_{j}) \in \Omega^{4j -1}(M)$.  
 
 \begin{thm}[Lazarov \cite{LazarovPasternack1976a, Lazarov1979}] \label{thm-gauge}
 Suppose that two framings $s, s'$ of $Q$ are related by a gauge transformation $\varphi \colon M \to \SOq$,   
$s'  = s \cdot \varphi$. 
Then on the level of forms,
\begin{equation}\label{eq-permanance1}
\Delta_{\F, g}^{s'}(h_{j}) = \Delta_{\F, g}^{s}(h_{j}) + \varphi^*(\tau_{j})
\end{equation}
In particular, for $j > q/4$, 
\begin{equation}\label{eq-permanance2}
\Delta_{\F, g}^{s'}[h_{j}] = \Delta_{\F,g}^{s}[h_{j}] + \varphi^*[\tau_{j}] \in H^{4j -1}(M ; \mR)
\end{equation}
 \end{thm}
 The relation (\ref{eq-permanance1}) can be used to easily calculate exactly how the cohomology classes 
  $\Delta_{\F, g}^{s}[h_I \otimes p_J]$ and $\Delta_{\F, g}^{s'}[h_I \otimes p_J]$  associated to   framings $s, s'$ are related.
  (See \S4 of \cite{LazarovPasternack1976a}, and \cite{Lazarov1979} for details.) 
  Here is one simple application of Theorem~\ref{thm-gauge}:
  
  {\bf Proof of Theorem~\ref{thm-loops}:} 
 For the product foliation of 
  $\SOq \times \mR^q$ we have a natural identification of the transverse orthogonal frame bundle $P = \SOq \times \SOq$. Let $s \colon \SOq \to P$ be the map $s(x) = x \times \{Id\}$, called  the product framing. Then   the   map $\ds \Delta_{\F,g}^s \colon RW_q   \to   \Omega^*(M) $ is identically zero. 
  
  On the other hand, the connecting map  $\delta \colon  \SOq \to \FRGq$ in (\ref{eq-puppe})  
classifies the   Riemannian foliation $\F_{\delta}$  of $\SOq \times \mR^q$, 
obtained via the pull-back of the standard product foliation of $\SOq \times \mR^q$  via the action  of    $\SOq$ on $\mR^q$.
 However, the normal  framing of  $\F_{\delta}$  is the product framing on 
$\SOq \times \mR^q$. Let $\varphi \colon \SOq \to \SOq$ be defined by $\varphi(x) = x^{-1}$ for $x \in \SOq$. 
Then $\F_{\delta}$ is diffeomorphic to the product foliation of $\SOq \times \mR^q$ 
with the framing defined by the gauge action of $\varphi$.

It follows from Theorem~\ref{thm-gauge} that for  $j > q/4$, 
$\Delta_{\F, g}^{s'}[h_{j}] =   \varphi^*[\tau_{j}]  = \pm {\tau_j} \in H^{4j -1}(\SOq ; \mR)$ is a generator. Hence, the connecting map $\delta \colon  \SOq \to \FRGq$ cannot be homotopic to the identity if there exists $j > q/4$ such that  ${\tau_j} \in H^{4j -1}(\SOq ; \mR)$ is non-zero. This is the case for all $q > 2$. \hfill $\eop$

The original Chern-Simons invariants of 3-manifolds \cite{ChernSimons1974} can be considered as examples of the above constructions.  Let $M$ be a closed oriented, connected 3-manifold with Riemannian metric $g$. Consider $M$ as foliated by points, then we obtain a Riemannian foliation of codimension $3$. Choose an oriented framing   $s \colon M \times \mR^3 \to TM$, then the transgression form $\Delta_{\F, g}^s (h_1) \in H^3(M ; \mR) \cong \mR$ is well-defined. Note that by formula (\ref{eq-permanance2}), the mod $\mZ$-reduction $\overline{\Delta_{\F, g}^s (h_1)} \in H^3(M ; \mR/ \mZ) \cong \mR/ \mZ$ is then independent of the choice of framing. This invariant of the metric is just the Chern-Simons invariant for $(M, g)$, as introduced in \cite{ChernSimons1974}.  On the other hand, Atiyah showed in  \cite{Atiyah1990} that for a 3-manifold, there is a ``canonical'' choice of framing $s_0$ for $TM$, so that there is a canonical $\mR$-valued Chern-Simons invariant,  $\Delta_{\F, g}^{s_0} (h_1) \in \mR$.

The paper of Chern and Simons also shows that the values of  
$\overline{\Delta_{\F, g}^s (h_1)} \in  \mR/ \mZ$ can vary non-trivially with the choice of Riemannian metric. Variational properties of the secondary classes are  discussed below in \S\ref{sec-variation}.

 One of the standard problems in foliation theory, is to determine whether the universal characteristic map is injective. For the classifying  space $B\Gamma_q$ of smooth foliations, this remains one of the outstanding open problems \cite{Hurder2008a}. In contrast, for Riemannian foliations,   the universal map (\ref{eq-universal}) is injective. We present here a new proof of this, based on Theorem~\ref{thm-thom2}, following the same idea of proof as in \cite{Hurder1981a}, but more explicit.

\begin{thm}[Hurder \cite{Hurder1981a}]  \label{thm-injectall}
There exists a compact  manifold  $M$ and a Riemannian foliation $\F$ of $M$ with trivial normal bundle, such that $\F$ is defined by a fibration over a compact manifold of dimension $q$, and  the   characteristic map
$  \Delta^s_{\F, g}  \colon H^*(RW_q) \rightarrow H^*(M)$ is injective. 
Moreover, if $q$ is odd, then $M$ can be chosen to be   connected. 
\end{thm}
\sop
 Let $B$ be the compact, oriented Riemannian manifold defined  in the proof   Theorem~\ref{thm-thom2}.    
Let $M$ be the bundle of oriented orthonormal frames for $TB$.  The basepoint map $\pi \colon M \to B$ defines a fibration
$\SOq  \to M \to B$, 
whose fiber $L_x = \pi^{-1}(x)$ over $x \in B$ is the group   $\SOq$ of  oriented  orthonormal frames in $T_xB$.
 Let  $\F$ be the foliation defined by the fibration. The Riemannian metric on $B$ lifts to the transverse metric on the normal bundle $Q = \pi^*TB$. 
 The lifted bundle $Q$ has a canonical framing $s$, where for $b \in B$ and $A \in \SOq$ the framing of $Q_{x,A}$ is that defined by the matrix $A$.   

   The normal bundle restricted to $L_x$ is trivial, as it is just the constant lift of $T_xB$.  That is, $Q | L_x \cong \pi^* (T_xB) \cong L_x \times \mR^q$. 
   The basic connection $\nabla_g$ restricted to $Q | L_x$ is the connection associated to the product bundle $L_x \times \mR^q$. 
However, the canonical framing of $Q \to M$ restricted to  $Q | L_x$   is twisted by $\SOq$.    Thus,  
the connection $\nabla_s$  on $Q$ for which  the canonical framing is parallel, restricts to  the Maurer-Cartan form on $\SOq \times \mR^q$ along each fiber $L_x$.

By Chern-Weil  theory, the forms
$\Delta_{\F,g}^s(h_{j}) =  h_{j}(\nabla_g,s)$   restricted to $L_x = \SOq$ are closed, and their classes in cohomology define the free exterior generators for the cohomology $H^*(\SOq ; \mR)$. (In the even case $q=2m$, one must include the Euler class $\chi_m$ as well.)

Give the algebra $RW_q$ the  basic filtration by the degree in $\cI(\SOq)_q$, and the forms in $\Omega^*(M)$ the basic filtration    by their   degree in $\pi^*\Omega^*(B)$. (See \cite{KT1975a} for example.)  Then the   characteristic map $\Delta_{\F,g}^s$ preserves the filtrations, hence induces a map of the associated Leray-Hirsch spectral sequences, 
$$ \Delta_r^{*,*} \colon E_r^{*,*}(RW_q, d_W)   \to   E_r^{*,*}(M, d_r) $$
For $r=2$,  we then have 
$$ \Delta_2^{*,*} \colon  E_2^{*,*}(RW_q) \cong (RW_q, d_W)  \to    E_r^{*,*}(M, d_2) \cong H^*(\SOq ; \mR) \otimes H^*(B ; \mR) $$
which is injective by the remark above. 
Pass to the $E_{\infty}$--limit to obtain that 
$$ \Delta_{\F, g}^s \colon H^*(RW_q)   \to  H^*(M) $$
induces an injective map of associated graded algebras, hence it is   injective. \hfill $\eop$

It seems to be an artifact of the proof that for $q \geq 4$ even, the manifold $M$ we obtain is not connected. 
\begin{prob}
For $q \geq 4$ even, does there exists a   closed, connected manifold  $M$ and a Riemannian foliation $\F$ of $M$ of codimension-$q$ and trivial normal bundle, such that   the secondary characteristic map  
$ \Delta_{\F, g}^s \colon H^*(RW_q)   \to  H^*(M) $ injects? Is there  a cohomological obstruction to the existence of such an example?
\end{prob}

Note that in the examples constructed in the proof of Theorem~\ref{thm-injectall}, the image of the     monomials  $h_I \otimes p_J$ for  $4i_1 + 4 |J|  > q$  (and hence $h_I \otimes p_J$ is $d_W$-closed) are integral:
$$\Delta_{\F, g}^s [ h_I \otimes p_J]   \in {\rm Image} \left\{H^*(M,\mZ) \to H^*(M,{\mR})\right\} $$
 This follows since the restriction of the forms $\Delta_{\F, g}^s(h_I)$ to the leaves of $\F$ are integral cohomology classes. In general, one cannot expect a similar integrality result to hold for examples with all leaves compact, as is shown by the Chern-Simons example previously mentioned. However,   a more restricted statement   holds.

\begin{defn}
A   foliation $\F$ of a   manifold $M$ is {\it compact Hausdorff} if every leaf of 
$\F$ is a compact manifold, and the leaf space $M/\F$ is a Hausdorff space.
\end{defn}

\begin{thm}[Epstein \cite{Epstein1976}, Millett \cite{Millett1975}] 
 A   compact Hausdorff foliation $\F$  admits a   holonomy-invariant Riemannian metric on its normal bundle $Q$. 
\end{thm}

In the next section, we discuss the division of the secondary classes into ``rigid'' and ``variable'' classes. One can show the following:
\begin{thm}
Let $\F$ be a compact Hausdorff foliation of codimension $q$  of $M$ with trivial normal bundle. 
If $h_I \otimes  p_J$ is a rigid class,  then  
$$\Delta_{\F, g}^s [h_I \otimes p_J]   \in {\rm Image}\left\{H^*(M,{\mQ}) \to H^*(M,\mR) \right\} $$
\end{thm}
This was proven in \cite{Hurder1980} for the case when  the leaf space $M/\F$ is a smooth manifold.  

 It is an interesting problem to determine geometric conditions on a Riemannian foliation which imply the rationality of the  secondary classes, as given for example by  Dupont and Kamber  in \cite{DK1993}. Rationality should be associated to rigidity properties for the global holonomy of the leaf closures \cite{Zimmer1988a}, one of the fundamental geometric concepts in the Molino Structure theory discussed in \S\ref{sec-molino}. One also expects rationality results for the secondary classes analogous to the celebrated results of Reznikov \cite{Reznikov1995, Reznikov1996}, possibly with some additional assumptions on the geometry of the leaves.

 \section{Variation of secondary   classes} \label{sec-variation}

The  secondary classes of a foliation are divided into two types, the  ``rigid''    and the ``variable'' classes.
Examples show that the variable classes are  sensitive to both the geometry and dynamical properties of the foliation,   while 
the rigid classes seem to be topological in nature.

A monomial $h_I \otimes  p_J \in RW_q$ is said to be \emph{rigid} if   $\deg  (p_{i_1} \wedge p_J) > q+2$.   Note that if $4i_1 + 4 |J|  > q$, then   this condition is automatically satisfied when $q = 4k$ or $q = 4k +1$. Here is the key property of the rigid classes:

\begin{thm}[Lazarov and Pasternack, Theorem~5.5 \cite{LazarovPasternack1976a}] \label{thm-rigid1}
Let $(\F_t , g_t, s_t )$ be a smooth 1-parameter family of framed Riemannian foliations.  Let  $h_I \otimes  p_J \in RW_q$ be a rigid class. Then
$$ \Delta^{s_0}_{\F_0 ,  g_0} [h_I \otimes p_J]  = \Delta^{s_1}_{\F_1 ,  g_1}  [h_I \otimes p_J]  \in H^*(M ; \mR)$$
Note that the family $\{(\F_t , g_t) \mid 0 \leq t \leq 1\}$ need not     be a  Riemannian foliation of codimension $q$ of $M \times [0,1]$.
\end{thm}

 For the special case where $q = 4k -2 \geq 6$, there is an extension of the above result:
 
 \begin{thm}[Lazarov and Pasternack, Theorem~5.6 \cite{LazarovPasternack1976a}] \label{thm-rigid2}
 Let $(\F, g_t, s_t)$ be a smooth 1-parameter family, where $\F$ is a fixed foliation of codimension $q$, each $g_t$ is a holonomy invariant    Riemannian metric on $Q$, and $s_t$ is a smooth family of framings on $Q$.  Let  $h_I \otimes  p_J \in RW_q$ satisfy $\deg  (p_{i_1} \wedge p_J) > q+1$.  Then
$$ \Delta^{s_0}_{\F ,  g_0} [h_I \otimes p_J]  = \Delta^{s_1}_{\F ,  g_1}  [h_I \otimes p_J]  \in H^*(M ; \mR)$$
We say that these classes are \emph{metric rigid}. Thus,   the classes $ [h_I \otimes p_J] \in H^*(RW_q)$  are metric rigid when $\deg  (h_{i_1} \otimes p_J) > q$, and   rigid under all deformations when $\deg  (h_{i_1} \otimes p_J) > q+1$.
\end{thm}

   A closed monomial $h_I \otimes p_J$ which is not rigid, is said to be \emph{variable}.  In the special case $q=2$, the class $[\chi_1 \otimes e_1] \in H^3(RW_2)$   is variable. For $q  > 2$,   neither the Euler class $e_m$ or its transgression $\chi_m$ can occur in a variable class, so for $q = 4k -2$ or $q = 4k -1$, the   variable classes  are spanned by the closed monomials
  \begin{equation}
\cV_{q}   =    \{ h_I \otimes p_J \mid 4   i_1 + 4   |J| = 4 k \} 
 \end{equation}
Let $v_q^k$ denote the dimension  of the subspace of $H^k(RW_q)$ spanned by the variable monomials. 

Theorem~\ref{thm-rigid2} implies that for codimension $q = 4k -2 \geq 6$, in order to continuously vary the value of a variable class $h_I \otimes p_J$ it is necessary to deform the underlying    foliation. For $q = 4k -1$, the value of variable class may (possibly)  be continuously varied by simply changing the transverse metric for the foliation, as seen in various  examples. We illustrate this with two examples.
 
\begin{ex}[Chern-Simons, Example 2 in \S6 of \cite{ChernSimons1974}]  \label{ex-CS}
\end{ex}
Consider $\mS^3$ as the Lie group ${\bf SU}(2)$ with Lie algebra spanned by 
$$ 
X = \left[\begin{array}{cc}i &0\\0 &-i\\ \end{array}\right] ,
Y =  \left[\begin{array}{cc}0 &i\\i &0\\ \end{array}\right] ,
Z =  \left[\begin{array}{cc}0 &-1\\1 &0\\ \end{array}\right]  $$
which gives   a framing $s$ of $T\mS^3$. 
   Let $g_u$ be the Riemannian metric on $\mS^3$ for which the parallel Lie vector fields $\{u \cdot X, Y, Z\}$ are an orthonormal basis. Let $\F$ denote  the point-foliation of $\mS^3$.
   Then $[h_1] \in H^3(RW_3)$ and for each $u > 0$,  we have $\ds \Delta^{s}_{\F ,  g_u} [h_1] \in H^3(\mS^3 ; \mR) \cong \mR$.
   
\begin{thm}[Theorem~6.9, \cite{ChernSimons1974}]
$\ds \frac{d}{du} |_{u=1} \left ( \Delta^{s}_{\F ,  g_u} [h_1] \right) \ne 0$
\end{thm}

One expects similar results also hold for other compact Lie groups    of dimension $4k -1 \geq 7$, although the author does not know of a published calculation of this.  

Chern and Simons also prove a fundamental fact about the conformal rigidity of the transgression classes, which carries over to Riemannian foliations as their calculations are all local.  

\begin{thm}[Theorem~4.5, \cite{ChernSimons1974}]\label{thm-confinv}
Let $(\F, g)$ be a Riemannian foliation of codimension $q = 4k -1$ of the closed manifold $M$. Let $s$ be a framing of the normal bundle $Q$. 
Let $\mu \colon M \to \mR$ be a smooth function, which is constant along the leaves of $\F$.  Define a conformal deformation  of  $g$ by setting $g_t = \exp(\mu(t)) \cdot g$. Then for all $[h_I \otimes p_J] \in H^*(RW_q, d_W)$ with $4 i_1 + 4|J| = q+1$,
$$ \Delta^{s}_{\F ,  g_t} [h_I \otimes p_J]  = \Delta^{s}_{\F ,  g}  [h_I \otimes p_J]  \in H^*(M ; \mR)$$
That is, the rigid secondary classes in codimension $q = 4k-1$ are conformal invariants.
\end{thm}
Combining Theorems~\ref{thm-rigid1}, \ref{thm-rigid2} and \ref{thm-confinv} we obtain:
\begin{cor}\label{cor-confinv}
The secondary classes of Riemannian foliations are conformal invariants.
\end{cor}

\begin{ex}[Lazarov-Pasternack, \cite{LazarovPasternack1976b}]  \label{ex-LP2}
\end{ex}
A modification of  the  original examples of Bott  \cite{Bott1967}  and Baum-Cheeger \cite{BaumCheeger1969}  show that all of the variable secondary classes vary independently, by a suitable variation of foliations.
 Let $\alpha = (\alpha_1, \ldots , \alpha_{2k}) \in \mR^{2k}$. Define a Killing vector field $X_{\alpha} $ on $\mR^{4k}$, with coordinates 
 $(x_1 , y_1 , x_2 , y_2 , \ldots , x_{2k}, y_{2k})$,
 $$X_{\alpha} = \sum_{i=1}^{k} ~ \alpha_i \{  x_i \partial / \partial y_i - y_i \partial /\partial x_i \} $$
 Let $\phi^{\alpha}_t \colon \mR^{4k} \to \mR^{4k}$ be the isometric  flow of $X_{\alpha}$, which restricts to an isometric flow on the unit sphere $\mS^{4k -1}$, so defines a Riemannian foliation $\F_{\alpha}$ of codimension $q = 4k -2$ of $\mS^{4k -1}$.
 
 Let $h_i \otimes p_J$ satisfy $4i + 4 |J| = 4k$. Associated to $p_i \wedge p_J$ is an Ad-invariant polynomial $\varphi_{i,J}$ on $\mathfrak{so}(4k)$ of degree $2k$. Let $M \to \mS^{4k -1}$ denote the bundle of orthonormal frames for the normal bundle to $\F_{\alpha}$, for $\alpha$ near $0 \in \mR^{2k}$. The spectral sequence for the fibration ${\bf SO}(4k -2) \to M \to \mS^{4k -1}$ collapses at the $E_2^{r,s}$-term, so we have an isomorphism $H^*(M ; \mR) \cong H^*(\mS^{4k -1} , \mR) \otimes H^*({\bf SO}(4k -2) ; \mR)$. Let $[C] \in H^{4k -1}(M, \mR)$ be the non-zero class corresponding to the fundamental class of the base. 
 
 \begin{thm}[\S\S2 \& 3, \cite{LazarovPasternack1976b}] 
 There exists $\lambda \ne 0$  independent of the choice of $p_i \wedge p_J$ such that  
 \begin{equation}
\langle \Delta^s_{\F_{\alpha}, g}[h_i \otimes p_J] , [C] \rangle =   \lambda \cdot  \frac{\varphi_{i,J}(\alpha_1, \ldots , \alpha_{2k})}{\alpha_1 \cdots \alpha_{2k} }
\end{equation}
 \end{thm}
 
 These examples are for $q = 4k -2$. Multiplying by a factor of $\mS^1$ in the transverse direction yields examples with codimension $4k -1$, and the same secondary invariants. Hence, we have the following corollary, due to Lazarov and Pasternack:
 
 \begin{cor} [Theorem~3.6, \cite{LazarovPasternack1976b}]
 Let $q = 4k -2$ or $4k -1$. Then evaluation on a basis of $H^{4k-1}(RW_q ; d_W)$ defines a surjective map
 \begin{equation}
\pi_{4k-1}(BR\Gamma^+_q) \to  {\bf R}^{v_q^{4k -1}} 
\end{equation}
In particular, all of the variable secondary classes in degree $4k -1$ vary independently.
 \end{cor}
 Although not stated in \cite{LazarovPasternack1976b}, these examples are sufficient to imply that all of the variable secondary classes for Riemannian foliations vary independently. This was stated as Theorem~4, \cite{Hurder1981b}.

 The reader may consult the papers 
 \cite{Hurder1981a, KT1975a, LazarovPasternack1976a, LazarovPasternack1976b, Mei1983, Yamato1979, Yamato1981} for a more extensive collection of examples of the calculation of the secondary classes for Riemannian foliations.   

\bigskip

\section{Molino Structure Theory} \label{sec-molino}

 In the previous section, it was observed that complete variation of the secondary classes can be obtained by deformations of the underlying Riemannian foliation. This suggests the problem of determining exactly what aspects of the dynamics of $\F$ contributes to the variation of the values of the secondary classes \cite{Molino1994}. 
The Molino Structure Theory  for Riemannian foliations provides a precise framework for studying this problem, as highlighted in Molino's  survey    \cite{Molino1994}. This theory  describes the dynamics and topology of a Riemannian foliation of a compact manifold. We recall below some of the main results; 
 the  reader can consult Molino \cite{Molino1982,Molino1988}, Haefliger \cite{Haefliger1985,Haefliger1989}, or  Moerdijk and Mr{\v{c}}un \cite{MM2003} for   further details.

Recall that we assume     $M$ is a closed, connected smooth manifold,   $(\F, g)$ is a smooth Riemannian foliation of codimension $q$  with  tangential distribution $F = T \F$, and that the normal bundle $Q \to M$ to $ \F$ is oriented.

 Let $\pi \, \colon \whM \to M$ denote the bundle of   orthonormal frames for $Q$ with positive orientation, where the fiber over $x \in M$ is  $\pi^{-1}(x) = \Fr^+(Q_x)$,   the space of orthogonal frames of $Q_x$ with positive orientation. The manifold  $\whM$ is a principal right  $\SOq$-bundle. Set  $\whx = (x,e) \in \whM$ for $e \in \Fr^+(Q_x)$.

The first  fundamental observation   is that there is a Riemannian foliation $\whF$ of $\whM$, whose leaves are the holonomy coverings of the leaves of $\F$, and such that $ \whF$ has no holonomy. The definition of $ \whF$ can be found in the sources cited above, but there is an easy intuitive definition.  Let $X$ denote a vector field on $M$ which is everywhere tangent to the leaves of $\F$, so that its flow $\varphi_t \colon M \to M$ defines  $ \F$-preserving diffeomorphisms. More precisely, for each $x \in M$, the path $t \mapsto \varphi_t(x)$ is a path in the leaf $\LH_x$ through $x$.  The     differential of these maps induce transverse isometries $D_x\varphi_t \colon Q_x \to Q_{\varphi_t(x)}$ which act on the oriented frames of $Q$, hence define paths in $\whM$. Given $\whx = (x,e) \in \whM$, the leaf $\whLH_{\whx}$ is defined by declaring that the path $t \mapsto D_x\varphi_t(e)$ is tangent to  $\whLH_{\whx}$. 
  It follows from the  construction that the restriction   $\pi \, \colon \whLH_{\whx} \to \LH_x$ of the projection $\pi$ to each leaf  of $\whF$ is a covering map.
  
There is an $\SOq$-invariant  Riemannian metric $\whg$ on $T\whM$ such that   $ \whF$ is   Riemannian. The metric $\whg$ satisfies     $d\pi \, \colon T \whF \to T \F$ is an isometry,  and the restriction of $\whg$  to the tangent space $T\pi$ of the fibers of $\pi$  is induced  from  the natural bi-invariant  metric on $\SOq$. Then 
   $d\pi$ restricted to the orthogonal complement $(T \whF \oplus T\pi)^{\perp}$ is a Riemannian submersion   to $Q$.

The second fundamental observation is that the foliation $ \whF$ of $\whM$ is \emph{Transversally Parallelizable} (TP). 
Let  ${\bf Diff(} \whM,  \whF{\bf )}$ denote  the subgroup of diffeomorphisms of $\whM$ which map leaves to leaves for $ \whF$, not necessarily taking a leaf to itself.   For example, given any vector field $\what{X}$ on $\whM$ which is everywhere tangent to the leaves of $ \whF$, then its flow $\what{\varphi}_t$ defines a 1-parameter subgroup of  
 ${\bf Diff(} \whM,  \whF{\bf )}$, which preserves the leaves themselves.
 The TP condition is that  ${\bf Diff(} \whM,  \whF{\bf )}$ acts transitively on $\whM$.

 Given $\whx = (x,e) \in \whM$, let $\oLHx$ denote the closure of the leaf $L_x$ in $M$, and let  $\owhLHx$ denote the closure    of the leaf $\whLH_{\whx}$ in $\whM$. For notational convenience, we set $N_x = \oLHx$ and $N_{\whx} = \owhLHx$.
 Note that the distinction between $N_x \subset M$ and  $N_{\whx} \subset \whM$    is indicated by the basepoint.
 \begin{prop} \label{prop-lciso}
 Given any pair of points $\whx, \why \in \whM$,  there is a diffeomorphism $\Phi \in {\bf Diff(} \whM,  \whF{\bf )}$ which restricts to a foliated diffeomorphism, $\Phi \colon N_{\whx} \to N_{\why}$.  Hence,  given any pair of points $x,y \in M$, the universal coverings of the leaves $\LH_x$ and $\LH_y$ of $ \F$ are diffeomorphic and  quasi-isometric.  
   \end{prop}
The     Molino structure theory gives   a   description of the closures of the leaves of $ \F$ and $ \whF$.
\begin{thm} [Molino \cite{Molino1982,Molino1988}]  \label{thm-molino1} 
Let  $ \F$ be a Riemannian foliation of a closed manifold $M$. 
\begin{enumerate}
\item For each $\whx \in \whM$, the leaf closure $N_{\whx}$ is a submanifold of $\whM$.
\item  The set of all   leaf closures  $N_{\whx} $ defines a foliation $\whE$ of $\whM$ with all leaves compact without holonomy.
\item  The quotient leaf space $\whW$ is a closed manifold with an induced right $\SOq$-action.
\item The  fibration  $\ds \whUp   \, \colon \whM \to \whW$, with fibers     the  leaves of $\whE$, is    $\SOq$--equivariant. 
\end{enumerate}
Let $\ds W = M/\overline{ \F}$ be  the quotient of $M$ by the closures of the leaves of $ \F$, and $\Up \, \colon M \to W$ the quotient map. 
\begin{enumerate}\setcounter{enumi}{4}
\item   \label{eq-square}
 $W$ is a Hausdorff space, and  there is an $\SOq$-equivariant commutative diagram:
\end{enumerate}
 \begin{equation*}
\begin{array}{rcccl}
&\SOq&{=}&\SOq&\\
 &\downarrow&&\downarrow& \\
&\whM&\stackrel{\whUp}{\longrightarrow}&\whW&\\
\pi&\downarrow&&\downarrow&\whpi\\
&M&\stackrel{\Up}{\longrightarrow}&W&
\end{array}
\end{equation*}
\end{thm}
The second result provides a description of the closures of the leaves of $ \F$ and $ \whF$, and  the structure of $ \whF | N_{\whx}$.
\begin{thm} [Molino \cite{Molino1982,Molino1988}] \label{thm-molino2}   
Let  $ \F$ be a Riemannian foliation of a closed manifold $M$. 
\begin{enumerate}
\item    
There exists a  simply connected  Lie group $G$, whose Lie algebra $\mathfrak{g}$ is spanned by the  holonomy-invariant vector fields on $N_{\whx}$ transverse to $\whF$, such that   the restricted foliation $ \whF$ of $N_{\whx}$ is a Lie $G$-foliation with all leaves dense,    defined by  a Maurer-Cartan  connection 1-form  $\ds \ds \omega^{\whx}_{\mathfrak{g}}  \, \colon T N_{\whx}  \longrightarrow  \mathfrak{g}$. \\
\item 
Let $\rho_{\whx} \colon \pi_1(N_{\whx}, \whx) \to G$ be the global holonomy map of the flat connection $\omega^{\whx}_{\mathfrak{g}}$. Then the image  $\whcNx \subset G$ of $\rho_{\whx}$ is dense in $G$. \\
 \end{enumerate}
\end{thm}

\section{Some open problems} \label{sec-problems}

Theorems~\ref{thm-molino1} and \ref{thm-molino2} suggest a number of questions about the secondary classes of Riemannian foliations.
It is worth recalling that for the example constructed in the proof of Theorem~\ref{thm-injectall} of a Riemannian foliation for which the characteristic map is injective,   all of its leaves are compact, and so the structural Lie group $G$ of Theorem~\ref{thm-molino2} reduces to the trivial group. For this example, all of the secondary classes are integral.
 
The first two problems invoke the structure of the quotient manifold $\whW = \whM/\cE$  and   space $\ds W = M/\overline{ \F}$. 
\begin{prob} \label{prob1}
Suppose that foliation $\overline{ \F}$ of $M$ by the leaf closures of $\F$ is a non-singular foliation. Show that all   secondary classes of $\F$ are rational. In the case where every leaf of $\F$ is dense in $M$, so $W$ reduces to a point,  what can be said about the values of the secondary classes?
\end{prob}
In all examples where there exists a family of foliations for which the secondary classes vary non-trivially, the quotient space $W$ is singular, hence the action of $\SOq$ on  $\whW$ has singular orbits. The action of $\SOq$ thus  defines a stratification of $\whW$. (See \cite{HurderToeben2008} for a discussion of the various stratifications.)
\begin{prob} \label{prob3}
How do the values of the secondary classes for a Riemannian foliation depend upon the $\SOq$-stratification of $\whW$? Are there conditions on the structure of the stratification which force the secondary classes to be rational?
\end{prob}
 The next   problems concern the role of the structural Lie group $G$ of a Riemannian foliation $\F$.
 
\begin{prob} \label{prob4}
Suppose the  structural Lie group $G$   is nilpotent. For example,  if all leaves of $\F$ have polynomial growth, the $G$ must be nilpotent \cite{Carriere1988,Zimmer1987}.  Show that all rigid secondary classes of $\F$ are rational. 
\end{prob}

All of the known examples of families of Riemannian foliations for which the secondary classes vary non-trivially  are obtained by the action of an abelian group $\mR^p$, and so the structural Lie group $G$ is necessarily abelian. In contrast, one can ask whether there is a generalization to the secondary classes of Riemannian foliations    of the results of Reznikov that the rigid secondary classes of flat bundles must be rational \cite{Reznikov1995,Reznikov1996}.
\begin{prob} \label{prob5}
Suppose   the  structural Lie group $G$   is semi-simple  with real rank at least $2$, without any factors of $\mR$.  Must the values of the secondary classes be rigid under   deformation? Are all of the characteristic classes of $\F$ are rational?
\end{prob}

\begin{prob} \label{prob6}
Assume the leaves of $\F$ admit a Riemannian metric for which they are Riemannian locally symmetric spaces of higher rank \cite{Zimmer1982, Zimmer1988a}. Must all of the characteristic classes of $\F$ be rational?
\end{prob}

The final question is more global in nature, as it asks how the topology of the ambient manifold $M$ influences the values of the secondary classes for a Riemannian foliation $\F$ of $M$. Of course, one influence might be that the cohomology group $H^{\ell}(M ; \mR) = \{0\}$ where $\ell = \deg(h_I \otimes p_J)$, and then $\Delta_{\F}(h_I \otimes p_J) = 0$ is rather immediate. Are there more subtle influences, such as whether particular restrictions on the fundamental group $\pi_1(M)$ restrict the values of the secondary classes for Riemannian foliations of $M$?

\begin{prob} \label{prob7} How does the topology of a compact manifold $M$ influence the secondary classes  
  for a Riemannian foliation $(\F, g)$ with  normal framing $s$ of $M$? 
\end{prob}

There are various partial results for Problem~\ref{prob7} in the literature \cite{LazarovPasternack1976a, Mei1983, Yamato1981}, but no systematic treatment.  It seems likely that an analysis such as in Ghys \cite{Ghys1984}  for Riemannian foliations of simply connected manifolds would yield new results in the direction of  this question.


\end{document}